\documentclass[reqno]{amsart}

\usepackage{graphicx}
\usepackage{mathrsfs}
\usepackage{mathtools}
\usepackage{amsmath}

\mathtoolsset{showonlyrefs=true}
\numberwithin{equation}{section}

\newtheorem{thm}{Theorem}[section]
\newtheorem{lem}[thm]{Lemma}
\newtheorem{pro}[thm]{Proposition}

\newtheorem{cor}[thm]{Corollary}

\newtheorem{claim}{Claim}

\def\R{{\mathbb R}}
\def\Z{{\mathbb Z}}
\def\pr{{\bf pr}}

\def\bb{\begin}
           \def\ea{\end{array}}
          \def\ec{\end{center}}
     \def\ed{\end{description}}
        \def\ee{\end{equation}}
       \def\eea{\end{eqnarray}}
\def\beaa{\bb{eqnarray*}}     \def\eeaa{\end{eqnarray*}}
 \def\et{\end{thebibliography}}

        \def\d{\delta}
      \def\e{\varepsilon}

\def\bar2{\doublebar}

\def\Sing{{\rm Sing}}

\def\dist{{\rm dist}}

\def\Rep{{\rm Rep}}

\def\Orb{{\rm Orb}}
\def\OS{{\rm OrientSh}}
\def\StSh{{\rm StSh}}

\def\Ps{{\rm Ps}}
\def\Pt{{\rm Pt}}

\begin{document}

\title{Oriented and standard shadowing properties for topological flows}

\author{Sogo Murakami}

\address{Graduate School of Mathematical Sciences, University of Tokyo, 3-8-1 Komaba, Tokyo, Japan}
\email{murakami-sogo880@g.ecc.u-tokyo.ac.jp}







\begin{abstract}
We prove that oriented and standard shadowing properties are equivalent for topological flows with finite singularites that are Lyapunov stable or backward Lyapunov stable.
Moreover, we prove that the direct product $\phi_1 \times \phi_2$ of two topological flows has the oriented shdowing property if $\phi_1$ with finite singuralities has the oriented shadowing property, while $\phi_2$ has the limit set consisting of finite singularities that are Lyapunov stable or backward Lyapunov stable.
\end{abstract}

\maketitle{}

\section{Introduction}
Shadowing has been an important notion in the qualitative  study of dynamical systems, which  goes back to the classical works in the hyperbolic theory (see \cite{S.P.} for instance) dealing with the standard shadowing property. 
Since then, relationships between several types of shadowing properties and the structural stability have been studied. 
Among them, Palmer, Pilyugin and Tikhomirov proved that the Lipschitz shadowing property and the structural stability are equivalent for $C^1$ vector fields on a compact manifold \cite{P.S.T.}.
Note that this  implies that every structurally stable vector field has the standard shadowing property. Moreover, 
Tikhomirov proved that the $C^1$ interior of the set of vector fields with the oriented shadowing property coincides with the set of structural stable systems on a  manifold of dimension at most $3$ \cite{Tikho3dim}. 

In this direction of finding sufficient conditions for having the standard shadowing property, 
it is worth mentioning a result  by  Komuro, who  proved that oriented and standard shadowing properties are equivalent for topological flows on a compact metric space  without singularities \cite{Komuro}. He also  posed a question whether these properties are still equivalent for topological  flows with singularities.
However,  in $2015$ Tikhomirov gave  an example of a vector field on a four-dimentional manifold having  the oriented shadowing property, but not having the standard shadowing property \cite{Tikho}.
In his  example, two hyperbolic singularities of saddle type have played an important role.

In this paper, we extend Komuro's result, proving  that oriented and standard shadowing properties are equivalent for topological flows with finite singulaities that are Lyapunov stable or  backward Lyapunov stable.
In addition, we prove that the direct product $\phi_1 \times \phi_2$ of two topological flows has the oriented shdowing property if $\phi_1$ with finite singuralities has the oriented shadowing property, while $\phi_2$ has the limit set consisting of finite singularities that are Lyapunov stable or backward Lyapunov stable.

Let $(M, \dist)$ be a compact metric space and
let $\phi$ be a topological flow on $M$ (i.e., $\phi : \R \times M \to M$ is a continuous map satisfying $\phi(0, x) = x$ and $\phi(s + t, x) = \phi(s, \phi(t, x))$ for all $s, t \in \R$ and $x \in M$).
For $x \in M$ and $r > 0$, define $B(p, x)$ to be
\[
B(p, x) = \{ y \in M ; \dist(x, y) < p \}.
\]
Denote by $\Sing(\phi)$ the set of singularities of $\phi$ (i.e., $\Sing(\phi) = \{ x \in M ; \phi(t, x) = x \text{ for all } t \in \R \}$).

We say that $p \in \Sing(\phi)$ is {\it Lyapunov stable} if for any neighborhood $V$ of $p$, there exists a neighborhood $U$ of $p$ such that
\[
\phi(t, x) \in V, \quad t \geq 0
\]
for all $x \in U$.
Similarly, $p \in \Sing(\phi)$ is said to be {\it backward Lyapunov stable} if for any neighborhood $V$ of $p$, there exists a neighborhood $U$ of $p$ such that
\[
\phi(t, x) \in V, \quad t \leq 0
\]
for all $x \in U$.
We say that $\xi  : \R \to M$ is a {\it $d$-pseudotrajectory} of $\phi$ if
\[
\dist ( \xi(t + s), \, \phi(s, \xi(t)) ) < d
\]
for all $t \in \R$ and $s \in [0, 1]$.
Let $\Ps_\phi(d)$ be the set of all $d$-pseudotrajectories of $\phi$.
When it is not necessary to specify $\phi$, we just write $\Ps(d)$ for simplicity.
Denote by $\Rep$ the set of all homeomorphisms from $\R$ to $\R$ which preserves the orientation.
For $\e > 0$, let
\[
\Rep(\e) = \{ f \in \Rep ; \lvert \frac{f(a) - f(b)}{a - b} - 1 \rvert < \e, \forall a, b \in \R, a > b \}.
\]

We say that a topological flow $\phi$ has the {\it standard shadowing property} if for every $\e > 0$ there exists $d > 0$ such that
if $\xi \in \Ps(d)$ then
\[
\dist ( \xi(t), \, \phi(h(t), x) ) < \e, \quad t \in \mathbb{R}
\]
for some $x \in M$ and $h \in \Rep(\e)$.
Denote the set of all topological flows with the standard shadowing property by $\StSh(M)$.
As a weaker form of shadowing properties, we say that a topoligocal flow $\phi$ has the {\it oriented shadowing property} if for every $\e > 0$ there exists $d > 0$ such that
if $\xi \in \Ps(d)$ then
\[
\dist ( \xi(t), \, \phi(h(t), x) ) < \e, \quad t \in \mathbb{R}
\]
for some $x \in M$ and $h \in \Rep$.
Denote the set of all topological flows with the oriented shadowing property by $\OS(M)$.

The main theorem is the following one:
\begin{thm}\label{thm.main}
If $\Sing(\phi)$ consists of finite number of Lyapunov stable or backward Lyapunov stable points,
then the following two conditions are equivalent:
\begin{itemize}
	\item $\phi \in \OS(M)$.
	\item $\phi \in \StSh(M)$.
\end{itemize}
\end{thm}
It is trivial that $\StSh(M) \subset \OS(M)$. For the proof of Theorem \ref{thm.main}, we need to show that $\OS(M) \subset \StSh(M)$. 
It is essential that the singuralities are either Lyapunov stable or backward Lyapunov stable because of the example by Tikhomirov mentioned above.
\paragraph{\it Remark}
We observe that there is an example of a topological flow on $B^2 = \{ x \in \R^2 ; \lvert x \rvert \leq 1 \}$ having the oriented shadowing property and a Lyapunov stable singularity.
The following construction of $f$ is close to that of $f$ in \cite[Section 2.6]{S.P.}.
Let $f_0 : [-1, 1] \to [-1, 1]$ be such that
\[
f_0(x) =
\begin{cases}
-1 + \frac{x + 1}{2}, & x \in [-1, -1/3],\\
2x, & x \in [-1/3, 1/3],\\
1 + \frac{x - 1}{2}, & x \in [1/3, 1]
\end{cases}
\]
and define $f : [0, 1] \to [0, 1]$ by $f(0) = 0$ and if $x \in (0, 1]$ then
\[
f(x) = 2^{-(n+2)} f_0 ( 2^{n+2}(x - 3 \cdot 2^{-(n+2)}) ) + 3 \cdot 2^{-(n+2)}, \quad x \in (2^{-(n+1)}, 2^{-n}]
\]
for some $n = 0, 1, \dotsb$.
Let $\varphi : \R \times [0, 1] \to [0, 1]$ be a topological flow on $[0, 1]$ such that $\varphi(1, x) = f(x)$.
Let $(r, \theta)$ be the polar coodinate where $r \in (0, 1]$ and $\theta \in (0, 2\pi]$.
For $x \in B^2 \setminus \{ (0, 0) \}$ with the polar coodinate $(r, \theta)$, let $\phi(t, x) \in B^2 \setminus \{ (0, 0) \}$ be a point with the polar coodinate 
\[
(\varphi(t, x), \theta + 2\pi t).
\]
Moreover, letting $\phi(t, (0, 0)) = (0, 0)$ for all $t \in \R$, we have a topological flow $\phi$ on $B^2$.
Then observe that for $n \geq 0$,
$\{ x \in \R^2 ; \lvert x \rvert = 2^{-n} \}$ is an attracting closed orbit and
$\{ x \in \R^2 ; \lvert x \rvert = 3 \cdot 2^{-(n + 2)} \}$ is a repelling closed orbit.
On the other hand, $(0, 0)$ is simultaneously Lyapunov stable and backward Lyapunov stable, but neither attacting nor repelling.
It is easy to check that this example has the oriented shadowing property,
and therefore, according to Theorem \ref{thm.main}, this flow has the standard shadowing property.

Let $(M_i, \dist_i)$ be a compact metric space, and denote  the limit set of $\phi_i$ by $L(\phi_i)$ for $i = 1, 2$ (for the definition of the limit set, see \cite[Definition 1.5.1]{F.H.}).
Then the second theorem can be stated as follows:
\begin{thm}\label{thm.2nd}
Let $\phi_1 \in \OS(M_1)$ and $\# \Sing(\phi_1) < \infty$.
Suppose that a topological flow $\phi_2$ on $M_2$ has $L(\phi_2)$ consisting of finite singularities that are Lyapunov stable or backward Lyapunov stable.
Then $\phi_1 \times \phi_2 \in \OS(M_1 \times M_2)$.
\end{thm}
In general, the direct product of two topological flows with the standard shadowing property does not necessarily have the oriented shadowing property.
For the example of such topological flow, see \cite[example 5]{Komuro}.

It is easy to see that if the direct product $\phi_1 \times \phi_2$ has the oriented shadowing property, then both $\phi_1$ and $\phi_2$ have the oriented shadowing property.
Thus we obtain the following corollary:
\begin{cor}
Let $\phi$ be a topological flow on a compact metric space $M$.
If $L(\phi)$ consists of finite singularities that are Lyapunov stable or backward Lyapunov stable, then $\phi \in \OS(M)$.
\end{cor}

In Section \ref{sec.pre}, a proposition and two lemmas are given.
The proposition is proved in Section \ref{sec.subthm} using one of the two lemmas.
We apply the proposition in order to prove Theorems \ref{thm.main} and \ref{thm.2nd} in Sections \ref{sec.mainthm} and \ref{sec.app}, respectively.

\section{Preliminaries}\label{sec.pre}
In this section, we suppose that $\phi$ is a topological flow on a compact metric space $M$.
Let $K \subset M \setminus \Sing(\phi)$ be a compact set.
Then we can take $T_0 \in (0, 1)$ such that
\begin{equation}
\phi(t, x) \neq x, \quad t \in [0, 2T_0] \label{eq.defT_0}
\end{equation}
for all $x \in K$ because $K \cap \Sing(\phi) = \emptyset$. 
Define 
\begin{equation}
\tilde{K} = \bigcap_{t \in [-2T_0, 2T_0]} \phi(t, K).\label{eq.tildeK}
\end{equation}
Denote by $\Pt$ the set of all functions $\xi : \R \to M$ satisfying
\[
\xi(t + nT_0) = 
\phi(t, \xi(nT_0))
\]
for all $t \in [0, T_0)$ and $n \in \Z$. 
For $d > 0$, let
\[
\Pt_\phi(d) = \{ \xi \in \Pt ; \dist ( \phi(T_0, \xi(nT_0) ), \, \xi((n + 1)T_0) ) < d, \quad n \in \Z \}.
\]
Similarly to $\Ps(d)$, we may write $\Pt(d)$ instead of $\Pt_{\phi}(d)$.
\\
\paragraph{\it Remark}
There exists a constant $L > 0$ such that if $\xi$ is a $d$-pseudotrajectory then
\[
\dist ( \xi(t), \, \xi'(t) ) < Ld, \quad t \in \R
\]
for some $\xi' \in \Pt(Ld)$.
Thus in order to prove that $\phi$ has the standard shadowing property, it is enough to show that for sufficiently small $\e_0 > 0$, there exists $d > 0$ such that for all $\xi \in \Pt(d)$ we have
\[
\dist ( \xi(t), \, \phi(h(t), x) ) < \e_0, \quad t \in \mathbb{R}
\]
for some $x \in M$ and $h \in \Rep(\e_0)$.

The following proposition claims that translating the oriented shadowing property to the standard property is possible along a finite part of a pseudotrajectory away from singularities, which is the main ingredient of the proofs of Theorems \ref{thm.main} and \ref{thm.2nd}.
\begin{pro}\label{pro.shadowing}
For every $\e > 0$, there exists $\e_1 > 0$ satisfying the following property:

If $y \in M$ satisfies
\[
\dist ( \xi(t + t_0), \, \phi(g(t), y) ) < \e_1, \quad t \in [0, T_1]
\]
and
\[
\xi(t + t_0), \, \phi(g(t), y) \in \tilde{K}, \quad t \in [0, T_1]
\]
for some $\xi \in \Pt$, $t_0 \in \R$, $T_1 \geq T_0$ and $g \in \Rep$ with $g(0) = 0$,
then there exists $\tilde{g} \in \Rep(\e)$ with $\tilde{g}(0) = g(0)$ and $\tilde{g}(T_1) = g(T_1)$ such that
\[
\dist ( \xi(t + t_0), \, \phi(\tilde{g}(t), y) ) < \e, \quad t \in [0, T_1].
\]
\end{pro}
\begin{lem}\label{lem.basic}
For $T_0 \in (0, 1)$ and $K$ given above, the following properties hold.
\begin{enumerate}
	\item[(a)] For every $\e \in (0, T_0)$, there exists $\delta > 0$ such that for all $x \in K$ and $s, t \in [0, T_0]$, if $\dist ( \phi(t, x), \, \phi(s, x) ) < \delta$, then $\lvert s - t \rvert < \e$.
	\item[(b)] For every $T \in (0, T_0)$, there exists $\rho > 0$ such that for all $x \in K$ and $t \geq 0$,
	if $\phi([0, t], x) \subset B(\rho, x)$, then $t \in [0, T]$.
\end{enumerate}
\end{lem}
Defining
\[
\hat{K} = \bigcap_{t \in [0, T_0]} \phi(-t, K), 
\]
we have the following lemma, which comes from Lemma \ref{lem.basic}.
\begin{lem}\label{lem.linearlize}
For every $\e \in (0, 1)$ and $T \in (0, T_0)$, there exists $\e' \in (0, \e)$ such that if
\[
x,y \in \hat{K}
\]
and $g \in \Rep$ with $g(0) = 0$ satisfy
\[
\dist ( \phi(t, x), \, \phi(g(t), y) ) \leq \e', \quad t \in [0, T),
\]
then $\lvert g(T)/T - 1 \rvert \leq \e$.
\end{lem}
The proofs of Lemmas \ref{lem.basic} and \ref{lem.linearlize} are essentially given in those of
\cite[Lemmas 3.3 and 3.4]{Komuro}, respectively.
In fact, the only difference of \cite[Lemmas 3.3 and 3.4]{Komuro}
from ours is that $x \in X$.
However, the $\phi$-invariance of $X$ is not essential in his proofs, whose arguments still work as well for $x \in K$ as long as $K$ is compact and away from singularities.
Lemma \ref{lem.linearlize} is proved from Lemma \ref{lem.basic}, but we use only Lemma \ref{lem.linearlize} for the proof of Proposition \ref{pro.shadowing} in the next section.

\section{Proof of Proposition \ref{pro.shadowing}}\label{sec.subthm}
In this section, we prove Proposition \ref{pro.shadowing}.
Fix $\e \in (0, 4/5)$.
There exists $T \in (0, T_0)$ with $T_0/T \in \Z$ such that
\begin{equation}
\dist ( x, \, \phi(t, x) ) < \frac{\e}{4}, \quad t \in [0, 2T] \label{eq.T}
\end{equation}
for all $x \in M$ and
\begin{align}
1 - \e 
\leq & (1 - \e/4)\frac{(1 - \e/4)(1 - 2T/T_0)}{1 + \e/4}\\
< & (1 + \e/4)\frac{(1 + \e/4)(1 + 2T/T_0)}{1 - \e/4} \leq 1 + \e. \label{eq.ineq}
\end{align}
By Lemma \ref{lem.linearlize},
there  exists $\e' \in (0, \e/4)$ such that if $x, y \in \hat{K}$
and $g \in \Rep$ with $g(0) = 0$ satisfy
\begin{equation}
\dist ( \phi(\tau, x), \, \phi(g(\tau), y) ) \leq \e', \quad \tau \in [0, T), 
\end{equation}
then
\begin{equation}
\lvert g(T)/T - 1 \rvert \leq \frac{\e}{4}. \label{eq.ee'}
\end{equation}
Since $\tilde{K}$ is compact, we can take $\e_1 > 0$ so that if $x, y \in \tilde{K}$ satisfy $\dist(x, y) \leq \e_1$, then
\begin{equation}
\dist (\phi(t, x), \, \phi(t, y) ) < \e', \quad t \in [-T, T]. \label{eq.e'}
\end{equation}

Now, as the hypothesis of Proposition \ref{pro.shadowing}, suppose $y \in M$ satisfies
\begin{equation}
\dist ( \xi(t + t_0), \, \phi(g(t), y) ) < \e_1, \quad t \in [0, T_1] \label{eq.hyp1}
\end{equation}
and
\begin{equation}
\xi(t + t_0), \, \phi(g(t), y) \in \tilde{K}, \quad t \in [0, T_1] \label{eq.hyp2}
\end{equation}
for some $\xi \in \Pt$, $t_0 \in \R$, $T_1 \geq T_0$ and $g \in \Rep$ with $g(0) = 0$.
Define $\hat{g} \in \Rep$ by
\[
\hat{g}(t) =
\begin{cases}
t, & t \in (-\infty, 0],\\
g(t), & t \in [0, T_1],\\
g(T_1) + t - T_1, & t \in [T_1, \infty).
\end{cases}
\]
Take $t_x, t_y \in [0, T)$ such that
\[
\frac{t_0 - t_x}{T}, \frac{t_0 + T_1 + t_y}{T} \in \Z,
\]
and let $t_0' = t_0 - t_x$ and $T_1' = T_1 + t_x + t_y$.
Set $N = T_1'/T$.
Then for $n = 0, 1, \dots, N-1$ and $\tau \in [0, T]$,
define
\[
g_n(\tau) = \hat{g}(\tau + t_0' + nT - t_0) - \hat{g}(t_0' + nT - t_0).
\]
Note that, by our choice of $T_1'$ and $t_0'$, we have $\xi(\tau + t_0' + nT) = \phi(\tau, \xi(t_0' + nT))$ for all $n = 0, \dots, N - 1$ and $\tau \in [0, T)$.
Then using \eqref{eq.e'}, \eqref{eq.hyp1} and \eqref{eq.hyp2}, we obtain
\begin{equation}
\dist ( \xi(\tau + t_0' + nT), \, \phi ( g_n(\tau), \phi( \hat{g}(t_0' + nT - t_0), y) ) ) < \e', \quad \tau \in [0, T) \label{eq.dd}
\end{equation}
and
\[
\xi(t_0' + nT), \, \phi(\hat{g}(t_0' + nT - t_0), y) \in \hat{K}.
\]
Thus \eqref{eq.ee'} implies
\begin{equation}
\lvert g_n(T)/T - 1 \rvert \leq \frac{\e}{4} \label{eq.g_n}
\end{equation}
for all $n = 0, \dots, N - 1$.
Define
\[
h(t_0' + nT - t_0) = \hat{g}(t_0' + nT - t_0)
\]
for each $0 \leq n \leq N-1$ and extend $h$ linearly between these points.
Then $h$ belongs to $\Rep(\e/4)$.
Let $t \in [t_0' + nT - t_0, t_0' + (n + 1)T - t_0]$.
Since
\[
h(t), \, \hat{g}(t) \in [\hat{g}(t_0' + nT - t_0), \, \hat{g}(t_0' + (n + 1)T - t_0)] \label{eq.set}
\]
and \eqref{eq.g_n},
we have
\begin{equation}
\lvert h(t) - \hat{g}(t) \rvert < 2T, \quad t \in [t_0' - t_0, t_0' + T_1' - t_0]. \label{eq.2T}
\end{equation}
Therefore, we can apply \eqref{eq.T} to have
\begin{equation}
\dist ( \phi(\hat{g}(t), y), \, \phi(h(t), y) ) < \e/4, \quad t \in [t_0' - t_0, t_0' + T_1' - t_0]. \label{eq.dist}
\end{equation}
From this together with \eqref{eq.dd}, it follows that 
\begin{align}
\dist ( &\xi(\tau + t_0' + nT), \, \phi(h(\tau + t_0' + nT - t_0), y) )\nonumber\\
&\leq \dist ( \xi(\tau + t_0' + nT), \, \phi(\hat{g}(\tau + t_0' + nT - t_0), y) )\nonumber\\
&+ \dist ( \phi(\hat{g}(\tau + t_0' + nT - t_0), y), \, \phi(h(\tau + t_0' + nT - t_0), y) )\nonumber\\
&\leq \e' + \e/4
\leq 2\e/4 \label{eq.xih}
\end{align}
for all $n = 0, 1, \dots, N - 1$ and $\tau \in [0, T)$.
Finally, define
$\tilde{g} : \R \to \R$ by
\[
\tilde{g}(t) = 
\begin{cases}
t, & t \leq 0,\\
\frac{g(T_1)}{h(T_1) - h(0)} (h(t) - h(0)), & t \in [0, T_1],\\
g(T_1) + t - T_1, & t \geq T_1.
\end{cases}
\]
Then, from \eqref{eq.g_n} and the fact that $h \in \Rep(\e/4)$, we have
\begin{align}
&\frac{(1 - \e/4)(1 - 2T/T_0)}{(1 + \e/4)}\\
&\leq \frac{(1 - \e/4)(1 - 2T/T_1)}{(1 + \e/4)}
\leq \frac{(1 - \e/4)(T_1 - 2T)}{(1 + \e/4)T_1}
\leq \frac{g(T_1)}{h(T_1) - h(0)}\\
&\leq \frac{(1 + \e/4)(T_1 + 2T)}{(1 - \e/4)T_1}
\leq \frac{(1 + \e/4)(1 + 2T/T_1)}{(1 - \e/4)}
\leq \frac{(1 + \e/4)(1 + 2T/T_0)}{(1 - \e/4)}.
\end{align}
Hence \eqref{eq.ineq} implies $\tilde{g} \in \Rep(\e)$.
Applying inequality \eqref{eq.2T}, we see that $\lvert h(0) - g(0) \rvert$, $\lvert h(T_1) - g(T_1) \rvert < 2T$.
Therefore, for $t \in [0, T_1]$,
\begin{align}
\lvert \tilde{g}(t) - h(t) \rvert
&= \Biggl\lvert \Biggl( \frac{g(T_1) - g(0)}{h(T_1) - h(0)} - 1 \Biggr) (h(t) - h(0)) + (g(0) - h(0)) \Biggr\rvert\nonumber\\
&= \Biggl\lvert \frac{g(T_1) - h(T_1) - (g(0) - h(0))}{h(T_1) - h(0)} (h(t) - h(0)) + (g(0) - h(0)) \Biggr\rvert\nonumber\\
&= \Biggl\lvert \frac{h(t) - h(0)}{h(T_1) - h(0)} (g(T_1) - h(T_1) - (g(0) - h(0))) + (g(0) - h(0)) \Biggr\rvert\nonumber\\
&= \Biggl\lvert \frac{h(t) - h(0)}{h(T_1) - h(0)} (g(T_1) - h(T_1))
+ \Biggl( 1 - \frac{h(t) - h(0)}{h(T_1) - h(0)} \Biggr) (g(0) - h(0)) \Biggr\rvert\nonumber\\
&\leq 2T. \label{eq.gh}
\end{align}
Thus, \eqref{eq.T}, \eqref{eq.xih} and \eqref{eq.gh} imply that
\beaa
\dist ( \xi(t + t_0), \, \phi(\tilde{g}(t), y) )
&\leq& \dist ( \xi(t + t_0), \, \phi(h(t), y) ) + \dist ( \phi(h(t), y), \, \phi(\tilde{g}(t), y) )\\
&<& \frac{2\e}{4} + \frac{\e}{4} < \e
\eeaa
for $t \in [0, T_1]$, finishing the proof.

\section{Proof of Theorem \ref{thm.main}.}\label{sec.mainthm}
In this section, we assume that $\phi \in \OS(M)$ and $\Sing(\phi)$ consists of finite number of Lyapunov stable or backward Lyapunov stable points.
For $p \in M$ and $\e > 0$, let
\begin{align}
\Omega^+(p, \e) &= \{ \xi(t) \in M ; \xi \in \Ps(\e), \, \xi(0) \in B(\e, p), \, t \geq 0 \},\\
\Omega^-(p, \e) &= \{ \xi(t) \in M ; \xi \in \Ps(\e), \, \xi(0) \in B(\e, p), \, t \leq 0 \}.
\end{align}
\begin{lem}\label{lem.stablesing}
If $p \in \Sing(\phi)$ is Lyapunov stable (resp. backward Lyapunov stable),
then
\[
\bigcap_{n \geq 1} \Omega^+(p, 1/n) = \{ p \} \quad
( \text{resp. } \bigcap_{n \geq 1} \Omega^-(p, 1/n) = \{ p \} )
.
\]
\end{lem}
\begin{proof}
If this does not hold for a backward Lyapunov stable point $p$, then we can choose
\[
x \in \bigcap_{n \geq 1} \Omega^+(p, 1/n)
\]
with $x \neq p$.
Let $r \in (0, \dist(p, x)/2)$.
Then, by the Lyapunov stability, there exists $\e \in (0, r)$ such that
\begin{equation}
\phi(t, x) \in B(r, p), \quad t \geq 0 \label{eq.stable}
\end{equation}
for all $x \in B(\e, p)$,
Take $n \geq 2/\e$ so that any $\xi \in \Ps(1/n)$ is $\e/2$-oriented shadowed.
Since $x \in \Omega^+(p, 1/n)$,
there exists $\xi \in \Ps(1/n)$ such that
$\xi(0) \in B(1/n, p)$ and
$\xi(t_0) = x$ for some $t_0 \geq 0$.
Since $\xi \in \Ps(\e/2)$, we have $h \in \Rep$ with $h(0) = 0$ and $p_0 \in M$ satisfying
\begin{equation}
\dist ( \xi(t), \, \phi(h(t), p_0) ) < \frac{\e}{2}, \quad t \in \R. \label{eq.lyap}
\end{equation}
Then \eqref{eq.lyap} and $\xi(0) \in B(1/n, p)$ imply that
\[
\dist (p_0, p) \leq \dist \big( p_0, \, \xi(0) ) + \dist ( \xi(0), \, p ) \leq \e/2 + 1/n \leq \e.
\]
Thus it follows from \eqref{eq.stable} that
\[
\phi(h(t_0), p_0) \in \Orb^+(B(\e, p)) \subset B(r, p).
\]
This inclusion and \eqref{eq.lyap} imply
\beaa
\dist(p, x)
&\leq& \dist ( p, \, \phi(h(t_0), p_0) ) + \dist ( \phi(h(t_0), \, p_0), x )\\
&=& \dist ( p, \, \phi(h(t_0), p_0) ) + \dist ( \phi(h(t_0), p_0), \, \xi(t_0) )\\
&<& r + \e/2 < 2r < \dist(p, x),
\eeaa
which is a contradiction.
The proof for a backward Lyapunov stable point is similar.
\end{proof}

Fix $\e_0  > 0$ such that
\begin{equation}
3\e_0 < \min \{ \dist(p, q) ; p, q \in \Sing(\phi), p \neq q \}. \label{eq.dfne_0}
\end{equation}
From Lemma \ref{lem.stablesing}, there exist $d_0 > 0$, $T_1 \in (0, 1)$ and a neighborhood $U_p$ of each $p \in \Sing(\phi)$ satisfying the following properties:
\begin{enumerate}
	\item[(i)] If $p$ is Lyapunov stable and $\xi \in \Ps(d_0)$ satisfies $\xi(0) \in U_p$, then $\xi(t) \in B(\e_0/4, p)$ for $t \geq -1$.
	\item[(ii)] If $p$ is backward Lyapunov stable and $\xi \in \Ps(d_0)$ satisfies $\xi(0) \in U_p$, then $\xi(t) \in B(\e_0/4, p)$ for $t \leq 1$.
	\item[(iii)] If $\xi(t_1) \in B(\e_0/4, p_1)$ and
	$\xi(t_2) \in B(\e_0/4, p_2)$ for some $\xi \in \Ps(d_0)$ and $p_1, p_2 \in \Sing(\phi)$ with $p_1 \neq p_2$, then $\lvert t_1 - t_2 \rvert \geq T_1$.
\end{enumerate}

Take $r_0 > 0$ satisfying $B(3r_0, p) \subset U_p$ for all $p \in \Sing(\phi)$, and set
\[
K = M \setminus \bigcup_{p \in \Sing(\phi)} B(r_0, p).
\]
Let $T_0 \in (0, T_1)$ be such that
\eqref{eq.defT_0} holds for all $x \in K$
and 
\begin{equation}
\tilde{K} \cup \bigcup_{p \in \Sing(\phi)} B(2r_0, p) = M, \label{eq.M}
\end{equation}
where $\tilde{K}$ is defined by \eqref{eq.tildeK}.

Apply Proposition \ref{pro.shadowing} with $\e$ replaced by $\e_0/4$ and let $\e_1 > 0$ be the constant given by the propositon for $\e = \e_0/4$.
Let $d_1 > 0$ be such that any $\xi \in \Pt(d_1)$ is $\min \{ r_0, \e_0/4, \e_1 \}$-oriented shadowed and $\Pt(d_1) \subset \Ps(d_0)$, and let $\xi \in \Pt(d_1)$.
Then the proof of Teorem \ref{thm.main} is reduced to showing that $\xi$ is $\e_0$-standard shadowed (see Remark in Section $2$).
The following claim is immediate.
\begin{claim}\label{claim.hasi}
Assume that $p$ is Lyapunov stable and $\xi(t) \in U_p$ hold for some $t \in \R$.
If $\inf \{ t \in \R ; \xi(t) \in U_p \} > -\infty$,
then there exists $t_0 \in \R$ such that
\begin{gather}
\xi(t) \in B(\e_0/4, p), \quad t \geq t_0 \text{ and }\\
\xi(t) \notin U_p, \quad t \leq t_0.
\end{gather}
\end{claim}
In fact, letting $t_0 = \inf \{ t \in \R ; \xi(t) \in U_p \} - 1/2$, this claim follows from property $({\rm i})$ above.
Define
\[
S(\xi) = \{ p \in \Sing(\phi) ; \xi(t) \in U_p \text{ for some } t \in \R \}. 
\]
Claim \ref{claim.hasi} shows that $S(\xi)$ contains at most one Lyapunov stable point
for otherwise we would have two Lyapunov stable points $p, q \in S(\xi)$
for which $\xi(t) \in B(\e_0/4, p) \cap B(\e_0/4, q)$ holds for sufficiently large $t \in \R$,
contradicting \eqref{eq.dfne_0}.
Similarly, $S(\xi)$ contains at most one backward Lyapunov stable point.

Let us consider several possible cases.
	\paragraph{\it Case 1} $S(\xi) = \emptyset$.

	There exist $h \in \Rep$ with $h(0) = 0$, $x_0 \in M$ and $t_0 \in \R$ such that
	\[
	\dist ( \xi(t + t_0), \, \phi(h(t), x_0) ) < \min \{r_0, \e_0/4, \e_1 \}, \quad t \in \R.
	\]
	Then it follows from \eqref{eq.M} that $\xi(t + t_0), \phi(h(t), x_0) \in \tilde{K}$ for all $t \in \R$.
	Then by applying Proposition \ref{pro.shadowing} for $[nT_0, (n + 1)T_0]$, $n \in \Z$ as $[0, T_0]$ in the proposition, there exists $\tilde{h} \in \Rep(\e_0)$ such that
	\[
	\dist ( \xi(t + t_0), \, \phi(\tilde{h}(t), x_0) ) < \e_0/4 < \e_0, \quad t \in \R.
	\]
	\paragraph{\it Case 2} $S(\xi) = \{ p \}$ with Lyapunov stable point $p$.

	If $\inf \{ t \in \R ; \xi(t) \in U_p \} = -\infty$,
	then $\xi(t) \in B(\e_0/4, p)$ for all $t \in \R$, and therefore
	\[
	\dist ( \xi(t), \, \phi(t, p) ) < \frac{\e_0}{4}, \quad t \in \R.
	\]

	On the other hand if $\inf \{ t \in \R ; \xi(t) \in U_p \} > -\infty$,
	then from Claim \ref{claim.hasi} we have $t_0 \in \R$ satisfying
	\begin{equation}
	\xi(t) \in B(\e_0/4, p), \quad t \geq t_0 \label{eq.t_00}
	\end{equation}
	and
	\begin{equation}
	\xi(t) \notin U_p, \quad t \leq t_0.\nonumber
	\end{equation}
	By the choice of $d_1 > 0$, there exist $h \in \Rep$ with $h(0) = 0$ and $x_0 \in M$ such that
	\begin{equation}
	\dist ( \xi(t + t_0), \, \phi(h(t), x_0) ) < \min \{r_0, \e_0/4, \e_1 \}, \quad t \in \R.\label{eq.min}
	\end{equation}
	Then it follows from \eqref{eq.M} that $\xi(t + t_0), \phi(h(t), x_0) \in \tilde{K}$ for all $t \leq 0$.
	Then applying Proposition \ref{pro.shadowing} again for $[-(n + 1)T_0, -nT_0]$, $n = 0, 1, \dotsb$, as $[0, T_1]$ in the proposition, we can find $\tilde{h} \in \Rep(\e_0)$ with $\tilde{h}(0) = 0$ such that
	\begin{equation}
	\dist ( \xi(t + t_0), \, \phi(\tilde{h}(t), x_0) ) < \frac{\e_0}{4}, \quad t \leq 0.\label{eq.til1}
	\end{equation}
	When $t \geq 0$, \eqref{eq.t_00} and \eqref{eq.min} show that
	\beaa
	\dist ( \xi(t + t_0), \, \phi(\tilde{h}(t), x_0) )
	&\leq& \dist ( \phi(\tilde{h}(t), \, x_0), \xi(h^{-1} \circ \tilde{h}(t) + t_0) )\\
	&+& \dist ( \xi(t + t_0), \, p )
	+ \dist ( p, \, \xi(h^{-1} \circ \tilde{h}(t) + t_0) )\\
	&<& \frac{\e_0}{4} + \frac{\e_0}{4} + \frac{\e_0}{4}.
	\eeaa
	This and \eqref{eq.til1} imply that
	\[
	\dist ( \xi(t + t_0), \, \phi(\tilde{h}(t), x_0) ) < \e_0, \quad t \in \R.
	\]
	\paragraph{\it Case 3} $S(\xi) = \{ p \}$ with backward Lyapunov stable point $p$.

	This case can be argued similarly to Case $2$.
	\paragraph{\it Case 4} $S(\xi) = \{ p, q \}$ with $p$ and $q$, Lyapunov stable and backward Lyapunov stable, respectively.

	Using Claim \ref{claim.hasi}, we have $t_0, t_1 \in \R$ such that
	\begin{gather}
	\xi(t) \in B(\e_0/4, p), \quad t \geq t_0, \\
	\xi(t) \notin U_p, \quad t \leq t_0
	\end{gather}
	and
	\begin{gather}
	\xi(t) \in B(\e_0/4, q), \quad t \leq t_1,  \\
	\xi(t) \notin U_q, \quad t \geq t_1.
	\end{gather}
	Then property $({\rm iii})$ implies $t_0 - t_1 \geq T_0$.
	By the choice of $d_1 > 0$, there exist $h \in \Rep$ with $h(0) = 0$ and $x_0 \in M$ such that
	\[
	\dist ( \xi(t + t_1), \, \phi(h(t), x_0) ) < \min \{r_0, \e_0/4, \e_1 \}, \quad t \in \R.
	\]
	Then it follows from \eqref{eq.M} that
	\[
	\xi(t + t_1), \, \phi(h(t), x_0) \in \tilde{K}, \quad t \in [0, t_0 - t_1].
	\]
	Now, by the choice of $\e_1$, we can apply Proposition \ref{pro.shadowing} to have $\tilde{h} \in \Rep(\e_0)$ with $\tilde{h}(0) = 0$ and $\tilde{h}(t_0 - t_1) = h(t_0 - t_1)$ such that
	\[
	\dist ( \xi(t + t_0), \, \phi(\tilde{h}(t), x_0) ) < \frac{\e_0}{4}, \quad t \in [0, t_0  -t_1].
	\]
	Finally, by the same argument as in Case $2$, we obtain
	\[
	\dist ( \xi(t + t_0), \, \phi(\tilde{h}(t), x_0) ) < \e_0, \quad t \in (-\infty, 0] \cup [t_0 - t_1, \infty).
	\]
	This completes the proof of Theorem \ref{thm.main}.

\section{Prooof of Theorem \ref{thm.2nd}}\label{sec.app}
Let $M = M_1 \times M_2$ and let $\pr_i : M \to M_i$ be the $i$-th projection map for $i = 1, 2$.
Define a distance function $\dist : M \times M \to \R$ by 
\[
\dist(x, y) = \max \{ \dist_1(\pr_1(x), \pr_1(y)), \dist_2(\pr_2(x), \pr_2(y)) \}.
\]
Let $\phi  = \phi_1 \times \phi_2$ with $\phi_i$, $i = 1, 2$ given in the hypothesis of Theorem \ref{thm.2nd}.
Note that if $\xi \in \Ps_\phi(d)$, then $\pr_i \circ \xi \in \Ps_{\phi_i}(d)$ for $i = 1, 2$ and $d > 0$.
Set $\xi_i = \pr_i \circ \xi$ for $i = 1, 2$.
Take $\e_0 > 0$ satisfying \eqref{eq.dfne_0}.
In the following lemma, we deal with a flow $\phi_2$ on $M_2$, claiming that $\xi_2 \in \Ps_{\phi_2}(d)$ with small $d$ stays near the singularities when $\lvert t \rvert$ is large enough.

\begin{lem}\label{lem.iii}
There exist $d_0 \in (0, \e_0/4)$, $T_1 \in (0, 1)$ and a neighborhood $U_p$ of each $p \in \Sing(\phi_2)$ satisfying properties $({\rm i})$, $({\rm ii})$ and $({\rm iii})$ in Section \ref{sec.mainthm} with $\phi$ replaced by $\phi_2$.
\end{lem}
\begin{proof}
Let $p \in \Sing(\phi_2)$ be Lyapunov stable.
Then there exists $\d > 0$ such that
\begin{equation}
\phi_2(t, x) \in B(\e_0/8, p), \quad t \geq 0 \label{eq.testesss}
\end{equation}
for every $x \in \overline{B(\d, p)}$.
Let $\d' > 0$ be such that
\[
\phi_2(t, x) \in B(\d/2, p), \quad t \geq 0
\]
for every $x \in B(\d', p)$.
By the choice of $\e_0$ and the fact that $\Sing(\phi_2) = L(\phi_2)$,
the $\omega$-limit set of $x$ is $\{ p \}$ for $x \in \overline{B(\d, p)}$.
Thus, there exist a neighborhood $U_x$ of $x$ in $\overline{B(\d, p)}$ and $T_x > 0$ such that
\[
\phi_2(T_x, y) \in B(\d', p)
\]
for all $y \in U_x$.
By the compactness of $\overline{B(\d, p)}$ we can take $T > 0$ such that
\begin{equation}
\phi_2(T, x) \in B(\d/2, p) \label{eq.TTT}
\end{equation}
for all $x \in \overline{B(\d, p)}$.
Let $d_0 \in (0, \e_0/4)$ be such that
\begin{equation}
\dist_2 ( \phi_2(t, \xi_2(t_0) ), \, \xi_2(t + t_0) ) < \min \{ \e_0/8, \d/2 \}, \quad t \in [0, T] \label{eq.texi}
\end{equation}
for every $\xi_2 \in \Ps_{\phi_2}(d_0)$ with $\xi_2(t_0) \in B(\d, p)$.
If $\xi_2(t_0) \in B(\d, p)$ hold for some $t_0 \in \R$,
then by \eqref{eq.texi} and \eqref{eq.testesss}, we have
\begin{align}
\dist_2 ( \xi_2(t + t_0), \, p )
&\leq \dist_2 ( \xi_2(t + t_0), \, \phi_2(t, \xi_2(t_0)) )
+ \dist_2 ( \phi_2(t, \xi_2(t_0)), \, p ) \\
&< \frac{\e_0}{8} + \frac{\e_0}{8} = \frac{\e_0}{4}
\end{align}
for $t \in [0, T]$.
Moreover, by \eqref{eq.texi} and \eqref{eq.TTT},
\begin{align}
\dist_2(\xi_2(T + t_0), p)
&\leq \dist_2(\xi_2(T + t_0), \phi_2(T, \xi_2(t_0))) + \dist_2(\phi_2(T, \xi_2(t_0)), p)\\
&< \frac{\d}{2} + \frac{\d}{2} = \d.
\end{align}
From these inequalities, we can make an inductive argument to have
\[
\dist_2(\xi_2(t + t_0), p) < \frac{\e_0}{4}, \quad t \geq 0.
\]
Let $\d'' \in (0, \d)$ be such that
\[
\xi_2(t) \in B(\e_0/4, p), \quad t \in [-1, 0]
\]
when $\xi_2 \in \Ps_{\phi_2}(d_0)$
with $\xi_2(0) \in B(\d'', p)$.
Then setting $U_p = B(\d'', p)$ we obtain property $({\rm i})$.

The proof of property $({\rm ii})$ is similar,
and that of property $({\rm iii})$ is trivial making $d_0 > 0$ smaller if necessary.
\end{proof}
Take $r_2 > 0$ so small that
\begin{equation}
B(4r_2, p) \subset U_p \label{eq.defr_2}
\end{equation}
for every $p \in \Sing(\phi_2)$.
Let 
\[
F = M_2 \setminus \bigcup_{p \in \Sing(\phi_2)} B(r_2, p).
\]

In the proof of the following claim, we need the hypothesis $\Sing(\phi_2) = L(\phi_2)$.

\begin{claim}
There exist $d_0' \in (0, d_0)$ and $S_0 > 0$ such that if $\xi_2 \in \Ps_{\phi_2}(d_0')$ with $\xi_2(s_0) \in F$ for some $s_0 \in \R$,
then we have a backward Lyapunov stable point $p$ and a Lyapunov stable point $q$ satisfying
\begin{gather}
\xi_2(t + s_0) \in B(r_2, p), \quad t \leq -S_0,\label{eq.Up}\\
\xi_2(t + s_0) \in B(r_2, q), \quad t \geq S_0  \label{eq.Uq}
\end{gather}
and
\begin{equation}
\dist_2 ( \xi_2(t + s_0), \, \phi_2(t, \xi_2(s_0)) ) < \min \{r_2, \e_0/4 \}, \quad t \in [-S_0, 2S_0]. \label{eq.V}
\end{equation}
\end{claim}
This claim can be shown as follows.
As proved in Lemma \ref{lem.iii},
there exist $r > 0$ and $d > 0$ satisfying the following properties:
\begin{enumerate}
	\item[(i')] If $q$ is Lyapunov stable and $\xi_2 \in \Ps_{\phi_2}(d)$ satisfies $\xi_2(0) \in B(r, q)$, then $\xi_2(t) \in B(r_2, q)$ for $t \geq 0$.
	\item[(ii')] If $p$ is backward Lyapunov stable and $\xi_2 \in \Ps_{\phi_2}(d)$ satisfies $\xi_2(0) \in B(r, p)$, then $\xi_2(t) \in B(r_2, p)$ for $t \leq 0$.
\end{enumerate}
Since $\Sing(\phi_2) = L(\phi_2)$, for every $x \in F$, the $\alpha$-limit set of $x$ is a backward Lyapunov stable point and the $\omega$-limit set of $x$ is a Lyapunov stable one.
Thus, there exists $S_0 > 0$ such that if $x \in F$ then
\begin{equation}
\phi_2([-S_0, 0), x) \cap B(r/2, p) \neq \emptyset \label{eq.forr}
\end{equation}
for some backward Lyapunov stable point $p$, and
\[
\phi_2((0, S_0], x) \cap B(r/2, q) \neq \emptyset
\]
for some Lyapunov stable point $q$.
Let $d_0' \in (0, d)$ be such that if $\xi_2 \in \Ps_{\phi_2}(d_0')$ with $\xi_2(s_0) \in F$ for some $s_0 \in \R$, then
\begin{equation}
\dist_2 ( \xi_2(t + s_0), \, \phi_2(t, \xi_2(s_0)) ) < \min \{ r_2, \e_0/4, r/2 \} \quad t \in [-S_0, 2S_0]. \label{eq.diffxi}
\end{equation}
By \eqref{eq.forr}, there exists $t_0 \in [-S_0, 0)$ such that
\[
\phi_2(t_0, \xi_2(s_0)) \in B(r/2, p)
\]
for some backward Lyapunov stable point $p$.
This and \eqref{eq.diffxi} imply that
\begin{align}
\dist_2 ( p, \, \xi_2(t_0 + s_0) )
&\leq \dist_2 ( p, \, \phi_2(t_0, \xi_2(s_0)) )
+ \dist_2 ( \phi_2(t_0, \xi_2(s_0)), \, \xi_2(t_0 + s_0) )\\
&< r/2 + r/2 = r.
\end{align}
From this, we have $\xi_2 (t_0 + s_0) \in B(r, p)$ and then the property ({\rm i'}) gives \eqref{eq.Up} because $t_0 \geq -S_0$.
The proof of \eqref{eq.Uq} is similar, and
\eqref{eq.V} follows from \eqref{eq.diffxi}.

In the following lemma, we deal with a flow $\phi_1 \in M_1$. Although there is no assumption on the singularities of $\phi_1$, the lemma claims that when $\xi_1 \in \Ps_{\phi_1}$ passes through a small neighborhood of a singularity,
$\xi_1(t)$ and its shadowing trajectory stay near the singularity for all $t$ in some interval of $\R$.
\begin{lem}\label{lem.M}
There exist $\e_1 > 0$, $d_1 \in (0, d_0')$ and a neighborhood $V_p$ of each $p \in \Sing(\phi_1)$ such that if $\xi_1 \in \Ps_{\phi_1}(d_1)$ satisfies $\xi_1(s_0) \in V_p$ for some $s_0 \in \R$ and
\[
\dist_1 ( \xi_1(t + s_0), \, \phi_1(h(t), x_0) ) < \e_1, \quad t \in \R
\]
for some $h \in \Rep$ with $h(0) = 0$ and $x_0 \in M_1$,
then there exist $t_0 < -S_0$ and $t_1 > S_0$ such that
\begin{equation}
\xi_1(t + s_0), \phi_1(h(t), x_0) \in B(\e_0/2, p), \quad t \in [t_0, t_1]  \label{eq.claim}
\end{equation}
and
\[
\min_{i = 0, 1} \lvert h(t_i) \rvert \geq S_0.
\]
\end{lem}
\begin{proof}
Take $r > 0$ and $d_1 \in (0, d_0')$ such that if $\xi_1 \in \Ps_{\phi_1}(d_1)$ satisfies $\xi_1(s_0) \in B(r, p)$
for some $s_0 \in \R$
then
\begin{equation}
\xi_1(t + s_0) \in B(\e_0/4, p), \quad t \in [-S_0 - 1, S_0 + 1]. \label{eq.s_0}
\end{equation}
Put $V_p = B(r/2, p)$ and $\e_1 = \min \{r/2, \e_0/4 \}$.
As in the hypothesis, suppose that $\xi_1(s_0) \in V_p$ 
and
\[
\dist_1 ( \xi_1(t + s_0), \, \phi_1(h(t), x_0) ) < \e_1, \quad t \in \R.
\]
Then $\xi_1(s_0), x_0 \in B(r, p)$.
Let
\[
t_1 = \max \{ S_0 + 1, \, h^{-1}(S_0) \}
\]
and let
\[
t_0 = \min \{ -S_0 - 1, \, h^{-1}(-S_0) \}.
\]
Then it follows from \eqref{eq.s_0} that either $\xi_1(t + s_0)$ or $\phi_1(h(t), x_0)$ is contained in $B(\e_0/4, p)$ for any $t \in [t_0, t_1]$.
Since $\e_1 \leq \e_0/4$, we have \eqref{eq.claim}.
In addition, $t_1 \geq h^{-1}(S_0)$ and $t_0 \leq h^{-1}(-S_0)$ yield the conclusion.
\end{proof}

Take $r_1 > 0$ satisfying
\[
B(3r_1, p) \subset V_p
\]
for every $p \in \Sing(\phi_1)$,
and let
\[
K = M_1 \setminus \bigcup_{p \in \Sing(\phi_1)} B(r_1, p).
\]
Then choose $T_0 \in (0, T_1)$ depending on $r_1$ so small that
\eqref{eq.defT_0} holds for all $x \in K$
and 
\begin{equation}
\tilde{K} \cup \bigcup_{p \in \Sing(\phi_1)} B(2r_1, p) = M_1,
\end{equation}
where $\tilde{K}$ is given by \eqref{eq.tildeK}.

Let $\tau_0 > 0$ be such that
\begin{equation}
\sup_{x \in M}\{\dist ( x, \, \phi(t, x) ) \} < \min\{ r_2, \e_0/4 \}, \quad t \in [0, \tau_0]. \label{eq.deftau_0}
\end{equation} 

Now apply Proposition \ref{pro.shadowing} to $(M_1, \phi_1)$ with $\e$ replaced by $\min\{r_1, r_2, \e_1, \e_0/4, \tau_0/(2S_0) \}$ and let
\[
\e_2 \in (0, \min\{r_1, r_2, \e_1, \e_0/4, \tau_0/(2S_0) \})
\]
be the constant given by the proposition for $\e = \min\{r_1, r_2, \e_1, \e_0/4, \tau_0/(2S_0) \}$.

Take $d_2 \in (0, \min \{d_1, r_2/2 \})$ satisfying $\Pt_{\phi}(d_2) \subset \Ps_{\phi}(d_1)$.
Making $d_2 > 0$ smaller if necessary, we may assume that every $\xi_1 \in \Pt_{\phi_1}(d_2)$ is $\e_2$-oriented shadowed and every $\xi_2 \in \Pt_{\phi_2}(d_2)$ satisfies
\begin{align}
\inf \{ t \in \R ; \xi_2(t) \in U_p \text{ for some } p \in \Sing(\phi_2) \} &= -\infty, \label{eq.backphi}\\
\sup \{ t \in \R ; \xi_2(t) \in U_p \text{ for some } p \in \Sing(\phi_2) \} &= \infty. \label{eq.forphi}
\end{align}
Let $\xi \in \Pt_\phi(d_2)$.
The proof of Teorem \ref{thm.2nd} is reduced to showing that $\xi$ is $\e_0$-oriented shadowed (see the Remark in Section $2$).

We need to consider two cases.
The first one is:
\paragraph{\it Case 1}
There exists $s_0 \in \R$ such that $\xi(s_0) \in B(3r_1, \Sing(\phi_1)) \times F$.

Let $p_1 \in \Sing(\phi_1)$ satisfy
\[
\xi(s_0) \in B(3r_1, p_1) \times F \subset V_{p_1} \times F.
\]
By the choice of $d_2 > 0$, we can take $h \in \Rep$ with $h(0) = 0$ and $x_0 \in M_1$ such that
\begin{equation}
\dist_1 ( \xi_1(t + s_0), \, \phi_1(h(t), x_0) ) < \e_2, \quad t \in \R. \label{eq.orient}
\end{equation}
Then Lemma \ref{lem.M} implies that there exist $t_0 < -S_0$ and $t_1 > S_0$ such that
\begin{equation}
\xi_1(t + s_0), \phi_1(h(t), x_0) \in B(\e_0/2, p_1), \quad t \in [t_0, t_1] \label{eq.claim22}
\end{equation}
and
\[
\min_{i = 0, 1} \lvert h(t_i) \rvert \geq S_0.
\]
Define $\hat{h} \in \Rep$ by
\[
\hat{h}(t) = 
\begin{cases}
h(t),	& t \leq t_0,\\
\frac{S_0 + h(t_0)}{S_0 + t_0}(t - t_0) + h(t_0),		& t \in [t_0, -S_0],\\
t, 		& t \in [-S_0, S_0],\\
\frac{h(t_1) - S_0}{t_1 - S_0}(t - S_0) + S_0,		& t \in [S_0, t_1],\\
h(t),	& t \geq t_1.
\end{cases}
\]

Then let us prove that
\begin{equation}
\dist ( \xi(t + s_0), \, \phi(\hat{h}(t), (x_0, \xi_2(s_0))) ) < \e_0, \quad t \in \R. \label{eq.concl1}
\end{equation}
By \eqref{eq.claim22}, we have
\begin{equation}
\dist_1 ( \xi_1(t + s_0), \, \phi_1(\hat{h}(t), x_0) )
\leq \dist_1 ( \xi_1(t + s_0), \, p_1 ) + \dist_1 ( p_1, \, \phi_1(\hat{h}(t), x_0) )
 < \e_0 \label{eq.naka}
\end{equation}
for $t \in [t_0, t_1]$.
From \eqref{eq.orient} and the choice of $\e_2$ it follows that
\[
\dist_1 ( \xi_1(t + s_0), \, \phi_1(\hat{h}(t), x_0) ) < \e_0, \quad t \in (-\infty, t_0] \cup [t_1, \infty).
\]
Combining these two inequalities, we get
\begin{equation}
\dist_1 ( \xi_1(t + s_0), \, \phi_1(\hat{h}(t), x_0) ) < \e_0, \quad t \in \R.\label{eq.1stxi}
\end{equation}

On the other hand, since $\xi_2(s_0) \in F$, there exist $p, q \in \Sing(\phi_2)$ such that
\eqref{eq.Up}, \eqref{eq.Uq} and \eqref{eq.V} hold.

Then, by \eqref{eq.Uq} and \eqref{eq.V},
\[
\xi_2(S_0 + s_0), \, \phi_2(h(S_0), \xi_2(s_0)) \in B(2r_2, q) \subset U_q,
\]
and therefore using Lemma \ref{lem.iii} with respect to property $({\rm i})$, we get
\begin{align}
\dist_2 ( \xi_2(t + s_0), \, \phi_2(\hat{h}(t), \xi_2(s_0)) ) 
&\leq \dist_2 ( \xi_2(t + s_0), \, q ) + \dist_2 ( q, \, \phi_2(\hat{h}(t), \xi_2(s_0)) ) \\
&< \frac{\e_0}{4} + \frac{\e_0}{4}
= \frac{2\e_0}{4}
\end{align}
for $t \geq S_0$.
Similarly using Lemma \ref{lem.iii} with respect to property $({\rm ii})$, we get
\[
\dist_2 ( \xi_2(t + s_0), \, \phi_2(\hat{h}(t), \xi_2(s_0)) ) < \frac{2\e_0}{4}, \quad t \leq -S_0.
\]
From these inequalities and \eqref{eq.V}, it follows that
\[
\dist ( \xi_2(t + s_0), \, \phi_2(\hat{h}(t), \xi_2(s_0)) ) < \e_0, \quad t \in \R,
\]
which together with \eqref{eq.1stxi} imply \eqref{eq.concl1}.

\paragraph{\it Case 2}
For every $t \in \R$, $\xi(t) \notin B(3r_1, \Sing(\phi_1)) \times F$.

Define
\[
S_2(\xi) = \{ p \in \Sing(\phi_2) ; \xi_2(t) \in U_p \text{ for some } t \in \R \}. 
\]
As proved in Section \ref{sec.mainthm}, $S_2(\xi)$ contains at most one Lyapunov stable point, and so does for backward Lyapunov stable point.
Therefore \eqref{eq.backphi} implies $S_2(\xi) \neq \emptyset$.

Since $\xi_1 \in \Pt_{\phi_1}(d_2)$, the choice of $d_2 > 0$ implies that there exist $g \in \Rep$ and $x_1 \in M_1$ such that
\[
\dist_1 ( \xi_1(t + s_0), \, \phi_1(g(t), x_1) ) < \e_2, \quad t \in \R.
\]
If $S_2(\xi) = \{ p \}$ for some $p \in \Sing(\phi_2)$, then using Lemma \ref{lem.iii} with respect to properties $({\rm i})$ and $({\rm ii})$, together with \eqref{eq.backphi} and \eqref{eq.forphi}, we obtain
\[
\xi_2(t) \in B(\e_0/4, p), \quad t \in \R.
\]
Hence
\[
\dist ( \xi(t + s_0), \, \phi(g(t), (x_1, p)) ) < \max \{ \e_2, \e_0/4 \} < \e_0, \quad t \in \R
\]
as required.
Thus, we may assume that $S_2(\xi) = \{ p, q \}$ holds for a backward Lyapunov stable point $p$ and a Lyapunov stable point $q$.

As proved in Claim \ref{claim.hasi}, we can take $s_0 \in \R$ satisfying
$\xi_2(s_0) \in B(\e_0/4, p) \setminus  U_p \subset F$.
By \eqref{eq.Up} and \eqref{eq.Uq}, there exist $-S_0 \leq t_1 < 0 < t_2 \leq S_0$ such that
\begin{gather}
\xi_2(t + s_0) \in F, \quad t \in [t_1, t_2],\\
\xi_2(t_1 + s_0) \in B(2r_2, p), \text{ and } \\ 
\xi_2(t_2 + s_0) \in B(2r_2, q). \label{eq.t_2side}
\end{gather}

Take $h \in \Rep$ with $h(0) = 0$ and $x_0 \in M_1$ satisfying
\begin{equation}
\dist_1 ( \xi_1(t + t_1 + s_0), \, \phi_1(h(t), x_0) ) < \e_2, \quad t \in \R. \label{eq.2xi1}
\end{equation}
Then we have
\[
\xi_1(t + t_1 + s_0), \, \phi_1(h(t), x_0) \in \tilde{K}, \quad t \in [0, t_2 - t_1]
\]
because $\e_2 < r_1$.
By Proposition \ref{pro.shadowing},
there exists $\tilde{h} \in \Rep(\min\{ r_1, r_2, \e_1, \e_0/4, \tau_0/(2S_0) \})$ with $\tilde{h}(0) = 0$ and $\tilde{h}(t_2 - t_1) = h(t_2 - t_1)$ such that
\begin{equation}
\dist_1 ( \xi_1(t + t_1 + s_0), \, \phi_1(\tilde{h}(t), x_0) ) < \min\{ r_1, r_2, \e_1, \e_0/4, \tau_0/(2S_0) \} \label{eq.2tildexi1}
\end{equation}
for $t \in [0, t_2 - t_1]$.
Now define $\hat{h} \in \Rep$ by
\[
\hat{h}(t) = 
\begin{cases}
h(t),			& t \leq 0,\\
\tilde{h}(t),	& t \in [0, t_2 - t_1],\\
h(t),			& t \geq t_2 - t_1.
\end{cases}
\]
Then the proof of Theorem \ref{thm.2nd} is reduced to showing the following inequality:
\begin{equation}
\dist ( \xi(t + t_1 + s_0), \, \phi(\hat{h}(t), (x_0, \xi_2(t_1 + s_0))) ) < \e_0, \quad t \in \R.  \label{eq.tescon}
\end{equation}
From \eqref{eq.2xi1} and \eqref{eq.2tildexi1}, we have
\[
\dist_1 ( \xi_1(t + t_1 + s_0), \, \phi_1(\hat{h}(t), x_0) ) < \e_0, \quad t \in \R.
\]
For $t \in [0, t_2 - t_1]$, since $\tilde{h} \in \Rep(\tau_0/(2S_0))$, we have
\[
\lvert t - \hat{h}(t) \rvert
\leq \frac{\tau_0}{2S_0}\lvert t \rvert
\leq \tau_0.
\]
Therefore, from \eqref{eq.V} and \eqref{eq.deftau_0} it follows that
\begin{align}
\dist_2 &( \xi_2(t + t_1 + s_0), \, \phi_2(\hat{h}(t), \xi_2(t_1 + s_0)) )\\
&\leq \dist_2 ( \xi_2(t + t_1 + s_0), \, \phi_2(t, \xi_2(t_1 + s_0)) )\\
&+ \dist_2 ( \phi_2(t, \xi_2(t_1 + s_0)), \, \phi_2(\hat{h}(t), \xi_2(t_1 + s_0)) )\\
&< \min\{ r_2, \e_0/4 \} + \min\{ r_2, \e_0/4 \}
= \min\{ 2r_2, \e_0/2 \} \label{eq.dfdf}
\end{align}
for $t \in [0, t_2 - t_1]$.
This and \eqref{eq.t_2side} imply that 
\begin{align}
\dist_2 &( \phi_2(\hat{h}(t_2 - t_1), \, \xi_2(t_1 + s_0)), q )\\
&\leq \dist_2 ( \phi_2(\hat{h}(t_2 - t_1), \xi_2(t_1 + s_0)), \, \xi_2(t_2 + s_0) )
+ \dist_2 ( \xi_2(t_2 + s_0), \, q ) \\
&< 2r_2 + 2r_2.
\end{align}
Thus, by \eqref{eq.defr_2},
\[
\xi_2(t_2 + s_0), \, \phi_2(\hat{h}(t_2 - t_1), \xi_2(t_1 + s_0)) \in U_p.
\]
Then, use Lemma \ref{lem.iii} with respect to property $({\rm i})$ to get
\begin{align}
\dist_2 &( \xi_2(t + t_1 + s_0), \, \phi_2(\hat{h}(t), \xi_2(t_1 + s_0)) )\\
&\leq \dist_2 ( \xi_2(t + t_1 + s_0), \, q )
+ \dist_2 ( q, \, \phi_2(\hat{h}(t), \xi_2(t_1 + s_0))) )\\
&\leq \frac{\e_0}{4} + \frac{\e_0}{4} \label{eq.tend}
\end{align}
for $t \geq t_2 - t_1$.
Similarly, by considering the backward direction, we obtain
\begin{align}
\dist_2 ( \xi_2(t + t_1 + s_0), \, \phi_2(\hat{h}(t), \xi_2(t_1 + s_0)) )
\leq \frac{\e_0}{2}, \quad t \leq 0.\label{eq.ttend}
\end{align}
Now \eqref{eq.dfdf}, \eqref{eq.tend} and \eqref{eq.ttend} imply \eqref{eq.tescon} as required, completing the proof of Theorem \ref{thm.2nd}.

\section*{acknowledgements}
I would like to express my appreciation to my advisor S. Hayashi for his constructive suggestions and continuous support.

\bibliographystyle{amsplain}

\end{document}